\documentclass[12pt,a4paper, notitlepage]{article}
\usepackage[a4paper] {geometry}
\usepackage{amsmath,amsfonts,amssymb,mathrsfs,amsthm}
\usepackage{multirow}
\usepackage[english]{babel}
\usepackage[all]{xy}
\usepackage{tikz}
\numberwithin{equation}{section}
\numberwithin{figure}{section}
\newtheorem{thm}{Theorem}[section]
\newtheorem{prop}[thm]{Proposition}

\theoremstyle{definition}
\newtheorem{defx}[thm]{Definition}

\newcommand{\R}{I\!\!R}

\newcommand{\h}{I\!\!H}
\DeclareMathSizes{12}{13}{10}{10}

\title{Cohomological Equation on the Discrete Heisenberg Group}
\author{ {\emph{M.E. Egwe}}
\\ murphy.egwe@ui.edu.ng
\\ Department of Mathematics, University of Ibadan, Ibadan, Nigeria}
\begin{document}

\maketitle

\begin{abstract}
Let $\h_1$ be the one-dimensional Heisenberg group. In this paper, we consider some aspects of discrete dynamical systems on $\h_1$ and give a condition for the solution of a cohomological equations on the group.
\end{abstract}

\section{Introduction}
Many problems which usually arise due to the consideration of certain forms of rigidity and stability of physical bodies are modeled using dynamical systems. A discrete dynamical system is a couple $(M, \gamma)$, where $M$ is a manifold and $\gamma$ a diffeomorphism of $M$. The dynamics is usually given by the diffeomorphism $\gamma$ on the manifold $M$. The most basic cohomological equation which usually arises is a first order linear difference equation of the form: \begin{equation} f - f \circ \gamma = g, \; \text{where} \; f,g \in C^\infty(M). \end{equation}
First, we require a structure on the Lie group $G$ in order to measure the effect of $\gamma$. Objective measurement requires that the structure be invariant on $G$ or be preserved by the action of any vector $a \in G$. In other words, the structure will need to be invariant under change of coordinates. So, we assume invariant Haar measure on the compact Lie group $G$. Invariant integration also follows on the invariant measures \cite{Rudin73}.

Let $\mathscr{H}_1$ be the 3-dimensional discrete Heisenberg group. This is considered as the $3 \times 3$ upper triangular integral matrices with diagonal entries given by
\[ \mathscr{H}_1 = \begin{pmatrix} 1 & x & z \\ 0 & 1 & y \\ 0 & 0 & 1 \end{pmatrix}, \; x, y, z \in \mathbb{Z} \]
By notation, this is written as the triple $(x,y,z)$ and the group law or operation is given by the matrix multiplication: \[(x', y', z')(x,y,z) = (x' + x, y' + y, z' + z + x'y). \]
The inverse element is then seen to be \[ (x,y,z)^{-1} =  (-x, -y -z + xy). \]
The central element of $\mathscr{H}_1$ is the set $\mathfrak{Z}=\{ (0, 0, z) : z \in \mathbb{Z} \}$. More basic facts about $\mathscr{H}_1$ can be seen in the following:
\begin{prop}
Let $x', y', z', x, y, z$ be any integers. The multiplication in $\mathscr{H}_1$ satisfies the following equations: \\
(a) $(x, y, z)^{-1} = (-x, -y, -z + xy)$;\\
(b) $(x', y', z')\cdot (x,y,z) \cdot (x,y,z)^{-1} = (x', y', z + y'x - xy)$; \\
(c) $[(x', y', z'), (x, y, z)] = (0, 0, y'z - z'y)$. \\
The centre $\mathfrak{Z}(\mathscr{H}_1)$ coincides with $0 \times 0 \times \mathbb{Z} \subseteq \R^3 = \mathscr{H}_1$ and $[\mathscr{H}_1, \mathscr{H}_1 ] = \mathfrak{Z}(\mathscr{H}_1)$.
\end{prop}
Therefore, $\mathscr{H}_1$ is nilpotent of class two, and the canonical exact sequence \[ [\mathscr{H}_1, \mathscr{H}_1] \rightarrowtail \mathscr{H}_1 \twoheadrightarrow \mathscr{H}_1 \] presents $\mathscr{H}_1$ as a central extension of $\mathbb{Z}^2$ by $\mathbb{Z}$. The group is generated by the three elements \[ g_1 = \begin{pmatrix} 1 & 0 & 0 \\ 0 & 1 & 1 \\ 0 & 0 & 1  \end{pmatrix}, g_2 = \begin{pmatrix} 1 & 1 & 0 \\ 0 & 1 & 0 \\ 0 & 0 & 1 \end{pmatrix}, \; \text{and} \;  g_3 = \begin{pmatrix} 1 & 0 & 1 \\ 0 & 1 & 0 \\ 0 & 0 & 1 \end{pmatrix}. \]
We have the relations: \[ \begin{pmatrix} 1 & x & z \\ 0 & 1 & y \\ 0 & 0 & 1 \end{pmatrix} = g_1^yg_2^xg_3^z = g_2^xg_1^yg_3^{z-xy} \; \text{and} \; yx = xyz. \]
An element $g_1^x g_2^y g_3^z$ commutes with $g_1$ (resp. with $g_2$) if and only if $z = 0$ (resp. $x = 0$). Therefore, the centre is generated by the element $z$.

From this, it is clear that the group $\mathscr{H}_1$ modulo its centre $\mathfrak{Z}(\mathscr{H}_1)$ is abelian. That is, $\mathscr{H}_1/\mathfrak{Z}(\mathscr{H}_1)$ is the reduced Heisenberg group. Other forms of the generalized Heisenberg group can be found in
\cite{Egwe1},\cite{BMEE1}\cite{Egwe4}\cite{EgweBA} and the references therein.

Let $M$ be a left $\mathscr{R}$-module. We define $X = \widehat{M}$ as the dual group of $M$ in the sense of harmonic analysis. For $\gamma \in \mathscr{H}_1$, we have a homomorphism $\sigma_\gamma$ on the compact group $X$ defined by \[ (\sigma_\gamma x, m) = (x, \gamma^{-1}m) \; \text{with} \; x \in X \; \text{and} \; m \in M, \] where $(\cdot, \cdot)$ denotes the dual pairing of the pair $X$ and $M$. This defines an (algebraic) action $\sigma$ of $\mathscr{H}_1$ on $X$.

\begin{defx}
We call a measurable set $B \subseteq X$ invariant if for all $\gamma \in \mathscr{H}_1$, the equation $\sigma_\gamma(B) = B$ is satisfied modulo null sets. The action is ergodic if all invariant subsets have either measure one or measure zero, with respect to the normalized \emph{Haar measure}.
The action is strongly mixing if for $f,g \in L^2$, the inner product satisfies $\langle f, g \circ \sigma_\gamma\rangle \to \langle f,1\rangle \langle 1,g\rangle$  if $\gamma$ tends to infinity in the topology of one-point compactification of $\mathscr{H}_1$, that means if it leaves all finite sets.
\end{defx} In what follows, we highlight the representations of $\mathscr{H}_1.$

\section{Representations}
Let $\mathfrak{Z}$ be the subgroup defined by \[ \mathfrak{Z} := \{(0, 0, z) : z \in \mathbb{Z} \} \] and let $\mathscr{S}$ be the normal subgroup given by \[\mathscr{S} = \{ (0,y,z) : y, z \in \mathbb{Z} \} \] so, $\mathscr{S} \simeq \mathbb{Z}^2$ and $G = \mathscr{S} \ltimes \Gamma$. If we take $\mathscr{S} = \mathbb{Z}^2$, and let $s = (n,k)$ and $g = (m', n', k')$, then there is a left action of $G$ on $\mathscr{S}$ given by multiplication and projection \[ g\cdot s = (n + n' + k + k' + m'n) \eqno{(1)} \]
Using this same action of $G$ on $\mathscr{S}$ as given in (1), define operators $U_g \in \mathcal{U}(\mathscr{H}_1)$ by \[ U_gf(s) = f(g^{-1}\cdot s) = f(n - n', k - k' - m(n-n')) \eqno{(2)} \] The map $U : G \to \mathcal{U}(\mathscr{H}_1)$ given by $g \mapsto U_g$ is a unitary representation of $G$.

The Fourier transform and its inverse are unitary maps between $\ell^2(\mathbb{Z})$ and $L^2(\mathbb{T}^2)$, where $\mathbb{T} = \{ z \in \mathbb{C} : |z| = 1 \}$, where \[ \mathcal{F}^{-1}f(z,w) = \underset{n,k \in \mathbb{Z}}\sum f(n,k) z^n w^k \] and \[ \mathcal{F} h(n,k) = \int_{\mathbb{T}^2} h(z,w)z^{-n}w^{-k} dzdw. \]
Conjugation of the operator $U$ by $\mathcal{F}^{-1}$ yields an equivalent representation of $G$ on $L^2(\mathbb{T}^2, \nu)$, where $\nu$ is the Lebesgue measure on $\mathbb{T}^2$. For $f \in L^2(\mathbb{T}),$  let
\begin{eqnarray*}
f: \mathcal{F}^{-1}U_g \mathcal{F}f(z,w) &=& \underset{n,k \in \mathbb{Z}}\sum U_gf(n,k)z^nw^k\\
&=& \underset{n,k \in \mathbb{Z}}\sum \int_{\mathbb{T}^2} f(z', w')z^{-(n-n')}w^{-(k-k'-m(k-k'))}dz'dw'z^nw^k\\
&=& f(zw^{m'}, w)z^{n'}w^{k'}.\\
 \end{eqnarray*}

The new representation elements act as multiplication operators in the $w$-variable, so $\mathcal{F}^{-1}U\mathcal{F}$ can be expressed as a direct integral of representations of $G$, each on space $L^2(\mathbb{T})$. Specifically, for each $w \in \mathbb{T}$, and $g = (m', n', k')$ in $G$, define a unitary operator $U^w_g$ on $L^2(\mathbb{T}^2)$ by \[ U^w_gf(z) = f(zw^{m'})z^{n'} w^{k'}. \]
The map $g \mapsto U_g^w$ is a unitary representation $U^w$ of $G$ on $L^2(\mathbb{T})$ and $\displaystyle \mathcal{F}^{-1}U\mathcal{F}  = \int^{\oplus}_{\mathbb{T}} U^w dw.$ (Cf. \cite{Kornelson}).

Consider the discrete Heisenberg group $\mathscr{H}$ as the semi-direct product of $\mathbb{Z}^2$ and $\mathbb{Z}$: \[ \mathscr{H} = \mathbb{Z}^2 \ltimes \mathbb{Z}; \; \alpha : \mathbb{Z} \to Aut(\mathbb{Z}^2).\] The $\mathscr{H}$ consists of the triples $((m,k),s)$ with the multiplication: \[ ((m,k),s) \star ((m',k'),s') = ((m,k) + \alpha^s(m', k')) = ((m + m', k + k' + ms), s + s') \] In particular, \[((m,k), 0) \star ((0,0), s) = ((m,k), s); \] The dual object for $\mathbb{Z}^2$ is the torus $\mathbb{T}^2$.

A pair $(\xi, \eta) \in \mathbb{T}^2 $ corresponds to the character $\chi : (m,k) \mapsto e^{2\pi i(mk + k\eta)}$. The torus is the $G$-space for the action $\chi h(m,k) = \chi(h \star ((m,k),0) \star h^{-1})$. The action of $((m,k), s)$ is defined by the formula \[ (\xi, \eta) \mapsto (\xi, \eta) \begin{pmatrix} 1 & 0 \\ s & 1 \end{pmatrix} = (\xi + s\eta, \eta). \]

Now using the Mackey machine \cite{Mackey} for general nilpotent groups, all the representations, we realize an induced representation on $\mathscr{H}$ through these characters on $L^2(\{ x_0, x_1, \cdots, x_{p-1} \})$ as \[ [\rho(h)f](x) = A(h,x) f(xh) = e^{2\pi i(m\xi + (k + jm)\eta + [\frac{s+j}{p}])} f(xh), \] \[ \rho((m,k),s)f](x_j) = e^{2\pi(m\xi + (k + jm)\eta + [\frac{s+j}{p}])}f(x(j+s)\mod p), j = \overline{0, p-1} \]

By choosing in $L^2(\mathscr{H}) = L^2(\{x_0, x_1, \cdots, x_{p-1})$ the base $\epsilon_o, \epsilon_1, \cdots, \epsilon_{p-1}$, where $\epsilon_j$ is the indicator of a point $x_j \in \mathscr{H}$. With respect to this base, our representation is defined by \[ \rho((m,k),s) : \epsilon_j \mapsto e^{2\pi i(m\xi + (k + jm)\eta + [\frac{s+j}{p}]\alpha)} \epsilon(j-s)\mod p; \; j= \overline{0, p-1} \eqno{(4)} \]

So all finite dimensional irreducible representations of $\mathscr{H}$  are of the form $(4)$ with some $\xi, \eta, \alpha \in [0,1)$ and the orbit $\chi = (\xi, \eta) \in \mathbb{T}^2$ consists of $p$ points. The action of $((m,k), s)$ on $\mathbb{T}^2$ is given by \[ (\xi, \eta) \mapsto (\xi + s\eta, \eta). \]

\section*{The Character of the Representation}
We shall consider the character of the representation which we shall need in the sequel as given in [ref].
For $s$ non divisible by $p$, one has $\chi_\rho((m,k),s) = 0$. So, \[ \chi_\rho((m,k), s) = \delta^o_{s \mod p} \sum \exp(2\pi i(m\xi + k\eta + jm\eta + [\frac{s+j}{p}]\alpha)).\]

For $s$ divisible by $p$, and $j \in \overline{0, p-1}$. Then $[\frac{s+j}{p}] = \frac{s}{p}$ and we have \[ \chi_\rho((m,k), s) = \begin{cases} p \cdot e^{2\pi i(m\xi + k\eta + \frac{s}{p}\alpha)} & \text{if} \; s=0\mod p \; \text{and} \; m = 0\mod p \\ 0 & \text{Otherwise} \end{cases} \]
For any automorphism $\phi$, defined by \[ \phi((m,k), s) = ((s+m, k + \frac{m(m-1)}{2} + sm), m),\] one has \[ \chi_{\rho\phi}((m,k),s) = \begin{cases} p \cdot \exp(2\pi i((s+m)\eta + (-k + \frac{m(m-1)}{2}) + sm)\eta + \frac{m}{p}\alpha) & \text{if}\; s \; \text{and} \; m = 0 \mod p \\ 0 & \text{Otherwise} \end{cases}. \]

 Next, the concept of discrete Sobolev space on the Heisenerg group will be in order here. A comprehensive treatment of this concept can be seen in [Princeton-Saenz-Marcinkiewicz Multipliers]
Let $f \in \ell^2(\mathbb{Z})$. We define the operator $\Delta$ as \[ \Delta f(k) = f(k) - f(k-1), k \in \mathbb{Z}. \] If $f$ has a sufficiently rapid decay at infinity, then \[ (\Delta)\widehat{f }(\xi) = \underset{k \in \mathbb{Z}}\sum (f(k) - f(k-1))e^{-2\pi ik\xi} \] \[ = (1 - e^{-2\pi i\xi})\hat{f}(k). \]
We can thus define the operator $(1 + |\Delta|)^\alpha, \alpha > 0$ via its Fourier transform by \[ ((1 + |\Delta|)^\alpha)\widehat{f } (\xi) = (1 + |1 - e^{-2\pi i\xi}|)^\alpha \hat{f}(\xi). \]

We define the \underline{discrete Sobolev space} $\ell^2_\alpha(\mathbb{Z})$ as the set of functions $f$ on $\mathbb{Z}$ such that $(1 + |\Delta|)^\alpha f \in \ell^2(\mathbb{Z})$, with norm \[\begin{array}{rcl}
 ||f||_{\ell^2_\alpha(\mathbb{Z})}&=& ||(1 + |\Delta|)^\alpha f||_{\ell^2(\mathbb{Z})}\\
 \ \\
 &=& \left(\displaystyle\int^1_0 \left|(1 + |1 - e^{-2\pi i\xi}|)^\alpha \hat{f}(\xi)\right|^2 d\xi \right)^{1/2} \end{array}\]
\begin{thm}[Princeton-Saenz-Marcinkiewicz Multipliers]
Let $f \in \ell^2_\alpha(\R)$, the standard Sobolev space of degree $\alpha$, for $\alpha > 1/2$, and set $g = f|_{\mathbb{Z}}$. Then $g\in\ell^2_\alpha(\mathbb{Z})$ and, moreover, for every $R \geq \varepsilon \geq 0$, \[ ||(1 + |R\Delta|)^\alpha g ||_{\ell^2} \leq C_\varepsilon \left(\int_{\R} \left| (1 + |2\pi\xi R|)^\alpha \hat{f}(\xi) \right|^2 d\xi \right)^{1/2}, \] where the constant $C_\varepsilon$ depends only on $\varepsilon.\;\;\;\;\;\;\;\;\;\Box$
\end{thm}
\section{The Cohomology of $\mathscr{H}_1$}
Let $\mathfrak{H}$ be the integer lattice on $2n+1$-dimensional Heisenberg group $N$ for fixed  $n \in \mathbb{Z}$. Then $\mathfrak{H}_n$ can be thought of as the set given by $\{ (x,y,z) : x,y \in \mathbb{Z}^n, z \in \mathbb{Z} \}$, with group law given by \[ (x,y,z)(x', y', z') = (x + x', y + y', z + z' + \langle x,y \rangle) \eqno{(4)} \] where $\langle \cdot, \cdot \rangle$ represents the inner product of $\R^n$ restricted to $\mathbb{Z}$. Thus, $\langle x,y \rangle = \overset{n}{\underset{i=1}\sum} x_iy_i $.
A matrix representation of $\mathfrak{H}$ is given by the embedding into $SL(n+2, \mathbb{Z})$: \[(x,y,z) = (x_1, \cdots, x_n, y_1, \cdots, y_n, z) \to \begin{pmatrix} 1 & x_1 & x_2 & \hdots & x_n & z \\ 0 & 1 & 0 & 0 & \hdots & y_1 \\ 0 & 0 & 1 & 0  & \hdots & y_2 \\ \hdots & \hdots& \hdots & \hdots & \hdots & \hdots \\ \hdots & \hdots & \hdots & \hdots & \hdots & y_n \\ 0 & \hdots & \hdots & \hdots & \hdots & 1 \end{pmatrix} \in SL(n+2, \mathbb{Z}). \]
It is well known that $\mathfrak{H}$ is co-compact in $N$ and the nilmanifold $N/\mathfrak{H}$ is a classifying space for $\mathfrak{H}$. The centre of $\mathfrak{H}$ is given by \[ \mathfrak{Z} = [\mathfrak{H}, \mathfrak{H}] = \{(0, 0, z) : z \in \mathbb{Z} \} \simeq \mathbb{Z} \] so that $\mathfrak{H}/\mathfrak{Z} = \mathbb{Z}^{2n}$.

Thus $\mathfrak{H}$ has he structure of a central extension \[ 0 \longrightarrow \mathbb{Z} \longrightarrow \mathfrak{H} \longrightarrow \mathbb{Z}^{2n} \longrightarrow 0. \] We obtain the cohomology of $\mathfrak{H}$ from $[**]$ as follows.
\begin{thm}\cite{LeePacker}
Let $n \in \mathbb{Z}^+$ and let $\mathfrak{H}$ be the Heisenberg group of order $2n+1$. Then the cohomology group of $\mathfrak{H}$ with coefficients in $\mathbb{Z}$ viewed as a trivial module are given by
\[\mathscr{H}^k(H,\mathbb{Z}) =
\begin{cases} \overset{[k/2]}{\underset{j-0}\bigoplus}(Z_j)^{\begin{pmatrix} 2n \\ k - 2j \end{pmatrix} - \begin{pmatrix} 2n \\ k - 2j - 2 \end{pmatrix}};& 0 \leq k \leq n \\
Z^{ \begin{pmatrix} 2n \\ n \end{pmatrix} - \begin{pmatrix} 2n \\ n - 2 \end{pmatrix}} \oplus \begin{Bmatrix} \overset{[(n+1)/2]}{\underset{j=1}\bigoplus}(Z_j)^{\begin{pmatrix} 2n \\ n + 1-2j \end{pmatrix} - \begin{pmatrix} 2n \\ n-1-2j \end{pmatrix}} \end{Bmatrix},&k = n+1 \\
Z^{ \begin{pmatrix} 2n \\ k-1 \end{pmatrix} - \begin{pmatrix} 2n \\ k+1 \end{pmatrix}} \oplus \begin{Bmatrix} \overset{[(2n-k+2)/2]}{\underset{j=1}\bigoplus}(Z_j)^{\begin{pmatrix} 2n \\ k + 2j-1 \end{pmatrix} - \begin{pmatrix} 2n \\ k + 2j \end{pmatrix}} \end{Bmatrix},& n + 2 \leq k \leq 2n+1 \\
0,& k\geq 2n + 2. \end{cases} \]
\end{thm}
The discrete Heisenberg fan is given by \[\Sigma = \Sigma^* \bigcup \{(0,\xi) \in \mathbb{Z}^2 : \xi \geq 0 \} \] where \[\Sigma^*  = \{(\lambda, \xi) \in \mathbb{Z}^2 : \lambda \neq 0, \xi = |\lambda|(2j + n), j \in \mathbb{Z}^+ \}. \]
Recall that $\Sigma$ is homeomorphic to the Gelfand spectrum. The functions are bounded in $\mathscr{H}_1$ if and only if $(\lambda, \xi)$ belongs to the Heisenberg fan [Ref].
\begin{defx}
Let $t = (t_1, t_2, \cdots, t_n)$ be a vector in $\R^n$ such that the subgroup generated by its projection on $\mathbb{T} = \R^n/\mathbb{Z}^n$ is dense in $\mathbb{T}^n$. (This implies in particular that the number $1, t_1, t_2, \cdots, t_n$ are linearly independent over $\mathbb{Q}$. \\
(a) We say that $t$ is \emph{Diophantine} if there exists real number $C,s\in \mathbb{R}^+$ such that \[ |1 - e^{2\pi i \langle k,t \rangle}| \geq \frac{C}{|k|^s} \; \text{for any } k \in \mathbb{Z}^n. \]
(b) We say $t$ is \emph{Liouville vector} if there exists $C,s\in \mathbb{R}^+$ with $k_s \in \mathbb{Z}^n$ satisfying \[ |1 - e^{2\pi i \langle k,t \rangle}| \leq \frac{C}{|k_s|^s}. \]
\end{defx}
In the case of $\mathscr{H}_1$, we see a Diophantine vector to be such that there exists $C, s\in \mathbb{R}^+$ with \[| 1 - \chi(t_p)| \leq \frac{C}{|p|^s}, \; \text{where} \; t_p = \langle t,p \rangle, t, p \in \mathbb{Z}^n, \] and $\chi(t)$ an irreducible unitary representation of $\mathscr{H}_1$ or the character of $\mathscr{H}_1$. The Liouville vector can be defined analogously.\\
\ \\
\begin{defx} Let $\mathfrak{S}_n$ be the $n$-periodic sequence of complex numbers and $h=\{h_i\}_{i=1}^\infty\in \mathfrak{S}_n.$ The discrete Fourier transform of $h$ is the sequence $(\mathfrak{F}\{h\})_k=\widehat{h}_k$\\ where $\widehat{h}_k=\displaystyle\sum_{i=1}^{n-1} h_i\overline{w}^{ik},$ where $w=\exp(2\pi i/n)$. The discrete Fourier series is defined by $y=\displaystyle\displaystyle\sum_{i=1}^{n-1}\widehat{h}_k\exp(2\pi ik/n).$
\end{defx}
Next, consider a linear functional \[ \varphi : C^\infty(H) \to \mathbb{C}\] defined by \[ \varphi(f)(u) = \underset{k \in \mathbb{Z}}\sum \Delta f_k(u) = \underset{k \in \mathbb{Z}}\sum f_k e^{2\pi i\langle k,u \rangle }\] and \[ \varphi(g)(u) = \underset{k \in \mathbb{Z}^n}\sum \Delta g_k(u) = \underset{k \in \mathbb{Z}^n}\sum g_k e^{2\pi i \langle k,u \rangle}. \]

Thus the reduced system of equation is \[ (1 - e^{2\pi i\langle k, u \rangle})f_k = g_k, \; k \in \mathbb{Z}^n \eqno{(*)} \]
The necessary condition for $\underset{k \in \mathbb{Z}^n}\sum \Delta g_k(u) = 0$ is that $g_k = 0$. Thus \[ f_k(u) = \begin{cases} 0 & \text{if} \; k = 0 \\ \frac{g_k}{1 - e^{2\pi i\langle k,u \rangle}} & \text{if} \; m \neq 0. \end{cases} \eqno{(1)} \]
This function is then formally given by its Fourier transform coefficients $(f_k)_{k \in \mathbb{Z}}$.\\
\textbf{Definition 2.1:} A function $f:\R^n\longrightarrow\R$ is said to be radial if there is a function $\phi$ defined on $[0,\infty)$ such that $f(x)=\phi(|x|)$ for almost every $x\in \R^n$.\\
Simple and classical examples of radial functions and their properties can be seen in for example \\cite{Egwe4}\cite{Egwe5}\cite{Egwe6}\cite{Egwe7}\\
Let $\rho$ be transformation on $\R^n$ and $x\in\R^n$.  Then $\rho$ is said to be orthogonal if it is a linear operator on $\R^n$ that preserves the inner product $\langle\rho x,\rho y\rangle = \langle x,y\rangle$ for all $x,y\in \R^n$.\\
A Schwartz function $\varphi$ is said to be radial if for all orthogonal transformations, $A\in O(n)$ (i.e., for all rotations on $\R^n$), we have
$$\varphi=\varphi\circ A.$$ We shall denote the set of all radial Schwartz functions by $\mathscr{S}_{rad}(\R^n).$

A distribution $u\in \mathscr{S}'(\R^n)$ is called radial if for all orthogonal transformations $A\in O(n),$ we have $$u=u\circ A.$$ This means
$$\langle u,\varphi\rangle=\langle u,\varphi\circ A\rangle$$ for all Schwartz functions $\varphi$ on $\R^n.$ We denote by $\mathscr{S}_{rad}'(\R^n)$ the space of all radial tempered distributions on $\R^n.$\\
\begin{thm}
Let $\gamma$ be the diffeomorphism of associated to a translation of $\mathscr{H}_1$ by the vector $a = (a_1, a_2, \cdots, a_n, t)$, where $a_1, a_2, \cdots, a_n $ are linearly independent over $\mathbb{Q}$. Suppose $a$ is either Diophantine or Liouville, then there exists at least one solution for the equation $f - f \circ \gamma = g$ and the space $\mathscr{D}_\gamma(\mathbb{T})$ of radial distribution has dimension $-2n+2$ and is generated by the Haar measure $dx = dx_1\otimes dx_2 \otimes \cdots \otimes dx_n \otimes dt$.
\end{thm}
\textbf{Proof:} Using the Fourier coefficients, the equation (1) i.e., $f-f\circ\gamma=g$ yields the system $(*).$ Then a necessary condition for $g$ to be of the form $f-f\circ\gamma$ is $g_0=0,$ which implies $I(g)=0.$ Suppose this is satisfied,then the solution is of the form (1).

Clearly, the Fourier coefficients define  $C^\infty$-functions. We thus need to verify that for $r,s\in \mathbb{N}$ and $i=1,\cdots,p$ we have $\|f\|_{r,s}^i<+\infty.$\\
To do this, let $\mathfrak{U}$ be the operator acting of the smooth function $f$ on $\R$ by
$$\mathfrak{U}f(\xi,\lambda)=2\displaystyle\int_0^1\partial_\xi f(\xi,2\lambda \mu,\lambda)(1-\mu)d\mu.$$
Let $f$ be in $\mathscr{D}(\R^2)$ with support in $\{(\xi,\lambda)\in\R^2:|\xi|\leq\epsilon\}.$ Then,
$$\|\mathfrak{U}f\|_{L^\infty(\R^2)}\leq 2\displaystyle\int_0^1\|\partial^2_\xi f\|_{L^\infty(\R^2)}(1-\mu)d\mu\leq C\|\partial^2_\xi f\|_{L^\infty(\R^2)}.$$
For any Liouville vector $a,$ we have that
\begin{eqnarray*}
\|\mathfrak{U}f\|_{L^\infty(\R^2)}&\leq & C_{|a|}\displaystyle\sum_{s+r\leq 2|a|-2}\|\partial^s_\lambda\partial^r_\xi f\|_{L^\infty(\R^2)}\\
\ \\
&\leq & C_{|a|}\displaystyle\sum_{s+r\leq 2|a|-2}\sum_{=0}^s\|\partial_\lambda^{s-k}\partial_\xi^{(r+2+k)}f\|_{L^\infty(\R^2)}\\
\ \\
&\leq & C_{|a|}\displaystyle\sum_{s+r<2|a|}\|\partial^s_\lambda\partial_\xi^rf\|_{L^\infty(\R^2)}<\infty.
\end{eqnarray*} \hfill{$\Box$}

\bibliographystyle{amsplain}

\providecommand{\bysame}{\leavevmode\hbox to3em{\hrulefill}\thinspace}

\end{document}